\documentclass[12pt]{article}
\usepackage{amsmath}
\usepackage{amssymb}
\usepackage{vatola}

\def\q{\quad}

\def\qtq#1{\q\t{#1}\q}

\def\t{\text}
\def\f{\frac}

\def\b{\binom}

\def\sls#1#2{(\f{#1}{#2})}

\def\Ls#1#2{\Big(\f{#1}{#2}\Big)}
\let \pro=\proclaim
\let \endpro=\endproclaim

\begin{document}

\par\q \centerline {\bf  A kind of orthogonal polynomials and related
 identities II}

 \par\q\newline \centerline{Zhi-Hong Sun}
\par\q\newline \centerline{School of Mathematical Sciences}
\centerline{Huaiyin Normal University}
 \centerline{Huaian, Jiangsu 223300, P.R. China}
  \centerline{Email: zhsun@hytc.edu.cn}
\centerline{URL: http://www.hytc.edu.cn/xsjl/szh} \par\q
\newline
{\bf Abstract}  For $n=0,1,2,\ldots$ let
$d_n^{(r)}(x)=\sum_{k=0}^n\binom{x+r+k}k\binom{x-r}{n-k}$. In this
paper we illustrate the connection between $\{d_n^{(r)}(x)\}$ and
Meixner  polynomials.  New formulas and recurrence relations for
$d_n^{(r)}(x)$  are obtained, and a new
 proof of the formula for $d_n^{(r)}(x)^2$ is also given. In
 addition, for $r>-\f 12$ and $n\ge 2$ we show that
 $d_n^{(r)}(x)>\frac{(2x+1)^n}{n!}>0$ for $x>-\frac 12$, and
 $(-1)^nd_n^{(r)}(x)>0$ for $x<-\f 12$.

\par\q
\newline {\bf Keywords}:   polynomials; three-term recurrence; hypergeometric series; identity
\newline{\bf MSC(2010)}: 33C05, 33C20, 33C45, 05A19
\section*{1. Introduction}
\par In [5] the author introduced the polynomials
$$d_n^{(r)}(x)=\sum_{k=0}^n\b{x+r+k}k\b{x-r}{n-k}\
(n=0,1,2,\ldots)\tag 1.1$$ as a generalization of Delannoy numbers
(The Delannoy number $D(n,m)$ is the number of lattice paths from
$(0, 0)$ to $(m,n)$, with jumps $(0, 1),\ (1, 1)$ or $(1, 0)$, and
$D(n,m)=d_n^{(0)}(m)$). Set
$D_n^{(r)}(x)=(-i)^nn!d_n^{(r)}(\f{ix-1}2)$. From [5] we know that
 $\{D_n^{(r)}(x)\}$ are orthogonal polynomials for $r>-\f 12$. By [5,
Theorem 2.1],
$$\sum_{n=0}^{\infty}d_n^{(r)}(x)t^n=\f{(1+t)^{x-r}}{(1-t)^{x+r+1}}\qtq{for}
|t|<1.\tag 1.2$$ By [5, Theorem 2.2], $\{d_n^{(r)}(x)\}$ satisfy a
three-term recurrence relation:
$$(n+1)d_{n+1}^{(r)}(x)=(1+2x)d_n^{(r)}(x)+(n+2r)d_{n-1}^{(r)}(x)\q (n\ge 0).\tag
1.3$$ From [5, Theorem 2.6] we know that for $r\notin\{-\f 12,-\f
22,-\f 32,\ldots\}$ and $n\in\{0,1,2,\ldots\}$,
$$d_n^{(r)}(x)^2=\b{n+2r}n\sum_{k=0}^n\f{\b{x-r}k\b{x+r+k}k\b{n+2r+k}{n-k}}
{\b{2r+k}k}4^k.\tag 1.4$$ In [5] the author also established the
following linearization formula:
$$d_m^{(r)}(x)d_n^{(r)}(x)=\sum_{k=0}^{\min\{m,n\}}\b{m+n-2k}{m-k}\b{2r+m+n-k}k
(-1)^kd_{m+n-2k}^{(r)}(x).\tag 1.5$$
\par In this paper we continue to investigate the properties of
$d_n^{(r)}(x)$. In Section 2 we prove that
$$d_n^{(r)}(x)=\sum_{k=0}^n\b{n+2r}{n-k}\b{x-r}k2^k\tag 1.6$$
and reveal the connection between $d_n^{(r)}(x)$ and Meixner
polynomials. In Section 3 we establish the two-term recurrence
formula for $d^{(r)}_n(x)$:
$$(n+1)d^{(r)}_{n+1}(x)=(x+r+n+1)d^{(r)}_n(x)+(-1)^n(x-r)d^{(r)}_n(-x).\tag 1.7$$
We also determine $d_n^{(r)}(\f 12),\ d_n^{(r)}(1),\ d_n^{(r)}(\f
32)$ and $d_n^{(r)}(2)$. In Section 4, we give a new
 proof of (1.4) and state that for
$a\not\in\{0,1,2,\ldots\}$,
$$\Big(\sum_{k=0}^n\f{\b nk\b xk}{\b ak}(-2)^k\Big)^2
 =\f {(-1)^n}{\b
 an}\sum_{k=0}^n\f{\b xk\b{a-x}k\b{n+k-a-1}{n-k}}{\b ak}4^k.\tag 1.8$$
 In the
case $a=-1$, (1.8) was given by Guo [3]. In Section 5, we prove that
$\f{d_n^{(r)}(x)d_{n-1}^{(r)}(x)}{1+2x}>0$ for $n\ge 2$, $x\not=-\f
12$ and $r>-\f 12$, and pose an interesting conjecture.

\section*{2. A new formula for $d_n^{(r)}(x)$ and connections with Jacobi and Meixner polynomials}
\par \q We first present a new formula for $d_n^{(r)}(x)$.
 \pro{Theorem 2.1}  For
$n=0,1,2,\ldots$ we have
$$d_n^{(r)}(x)=\sum_{k=0}^n\b{n+2r}{n-k}\b{x-r}k2^k.$$
\endpro
Proof.  For $0<t<\f 13$, using Newton's binomial theorem we see that
 $$\align
 &\sum_{n=0}^{\infty}\sum_{k=0}^n\b{n+2r}{n-k}\b{x-r}k2^kt^n
 \\&=\sum_{k=0}^{\infty}\b{x-r}k(2t)^k\sum_{n=k}^{\infty}\b{n-k+k+2r}{n-k}t^{n-k}
 \\&=\sum_{k=0}^{\infty}\b{x-r}k(2t)^k\sum_{s=0}^{\infty}\b{k+2r+s}st^s
 \\&=\sum_{k=0}^{\infty}\b{x-r}k(2t)^k\sum_{s=0}^{\infty}\b{-k-2r-1}s(-t)^s
\\&=\sum_{k=0}^{\infty}\b{x-r}k(2t)^k(1-t)^{-k-2r-1}
\\&=(1-t)^{-2r-1}\sum_{k=0}^{\infty}\b{x-r}k\Ls {2t}{1-t}^k
\\&=(1-t)^{-2r-1}\Big(1+\f{2t}{1-t}\Big)^{x-r}=\f{(1+t)^{x-r}}{(1-t)^{x+r+1}}.
\endalign$$
 This together with (1.2) yields the result.

\pro{Corollary 2.1} Let $r\not\in\{-\f 12,-\f 22,-\f 32,\ldots\}$.
For $n=0,1,2,\ldots$ we have
$$\sum_{k=0}^n\b nk\f{d_k^{(r)}(x)}{\b{-2r-1}k}
=\f{\b{x-r}n}{\b{-2r-1}n}2^n.$$
\endpro
Proof. Observe that $\b {-a}k=(-1)^k\b{a+k-1}k$ and $\b an\b nk=\b
ak\b{a-k}{n-k}$. We get $\b{-2r-1}n\b
nk=\b{-2r-1}k\b{-2r-1-k}{n-k}=(-1)^{n-k}\b{-2r-1}k\b{n+2r}{n-k}$.
Hence from Theorem 2.1 we see that
 for
$r\not\in\{-\f 12,-\f 22,-\f 32,\ldots\}$,
$$ \f{d_n^{(r)}(x)}{(-1)^n\b{-2r-1}n}=\sum_{k=0}^n\b nk(-1)^k
\f{\b{x-r}k}{\b{-2r-1}k}2^k.\tag 2.1$$  Now applying
  the
binomial inversion formula we deduce the result.

\pro{Corollary 2.2} Let $r\not\in\{-\f 12,-\f 22,-\f 32,\ldots\}$.
For $n=0,1,2,\ldots$ we have
$$\sum_{k=0}^{[n/2]}\b n{2k}\f{d_{2k}^{(r)}(x)}
{\b{-2r-1}{2k}}=\f{\b{x-r}n+\b{-1-x-r}n}{\b{-2r-1}n}2^{n-1}$$ and
$$\sum_{k=0}^{[(n-1)/2]}\b n{2k+1}\f{d_{2k+1}^{(r)}(x)}
{\b{-2r-1}{2k+1}}=\f{\b{x-r}n-\b{-1-x-r}n}{\b{-2r-1}n}2^{n-1},$$
where $[a]$ is the greatest integer not exceeding $a$.
\endpro
Proof.  By [5, (2.3)], $d_n^{(r)}(x)=(-1)^nd_n^{(r)}(-1-x).$ Thus,
using Corollary 2.1 we see that
$$\sum_{k=0}^n\b nk(-1)^k\f{d_k^{(r)}(x)}{\b{-2r-1}k}=
\sum_{k=0}^n\b nk\f{d_k^{(r)}(-1-x)}{\b{-2r-1}k}
=\f{\b{-1-x-r}n}{\b{-2r-1}n}2^n.\tag 2.2$$ Combining Corollary 2.1
with (2.2) yields the result.

\par Let $P^{(\alpha,\beta)}_n(x)$ be the Jacobi polynomial given by
$$P^{(\alpha,\beta)}_n(x)=\f
1{2^n}\sum_{k=0}^n\b{n+\alpha}k\b{n+\beta}{n-k}(x+1)^k(x-1)^{n-k}.$$
Then
$$P^{(x-r-n,2r)}_n(3)=\sum_{k=0}^n\b{x-r}k\b{n+2r}{n-k}2^k.$$
This together with Theorem 2.1 and the facts
$P^{(\alpha,\beta)}_n(-x)=(-1)^nP^{(\beta,\alpha)}_n(x)$ (see [1-2])
and $d^{(r)}_n(-1-x)=(-1)^nd^{(r)}_n(x)$ yields the following
result. \pro{Corollary 2.3} For $n=0,1,2,\ldots$ we have
$$d^{(r)}_n(x)=P^{(x-r-n,2r)}_n(3)=(-1)^nP^{(2r,x-r-n)}_n(-3)
=P^{(2r,-1-x-r-n)}_n(-3).$$
\endpro
\par  Let $(a)_0=1$ and $(a)_k=a(a+1)\cdots(a+k-1)$ for
$k=1,2,3,\ldots$. The Meixner polynomials $\{M_n(x;b,c)\}$ are
defined by (see [1-2])
$$M_n(x;b,c)=\sum_{k=0}^n\f{(-n)_k(-x)_k}{(b)_k\cdot
k!}\Big(1-\f 1c\Big)^k.$$ Note that $(a)_k=(-1)^k\b {-a}kk!=\b
{a+k-1}kk!$. From (2.1) we deduce the following result.
\pro{Corollary 2.4} Let $r\not\in\{-\f 12,-\f 22,-\f 32,\ldots\}$.
For $n=0,1,2,\ldots$ we have
$$d^{(r)}_n(x)=\f{(2r+1)_n}{n!}M_n(x-r;2r+1,-1).$$
\endpro

\section*{3. New recurrence relations for $d^{(r)}_n(x)$ and special values of $d^{(r)}_n(x)$}

\par  We first present a
two-term recurrence formula for $d^{(r)}_n(x)$, which is better than
the three-term recurrence relation (1.3).
\pro{Theorem 3.1} For
$n=0,1,2,\ldots$ we have
$$(n+1)d^{(r)}_{n+1}(x)=(x+r+n+1)d^{(r)}_n(x)+(-1)^n(x-r)d^{(r)}_n(-x).$$
\endpro
Proof.  Let $f'(t)$ be the derivative of $f(t)$. By (1.2), for
$0<t<1$ we have
$$\align
&\sum_{n=0}^{\infty}\big((n+1)d^{(r)}_{n+1}(x)-nd_n^{(r)}(x)-(x+r+1)d_n^{(r)}(x)-(-1)^n(x-r)d^{(r)}_n(-x)\big)t^n
\\&=\Big(\sum_{n=0}^{\infty}d^{(r)}_{n+1}(x)t^{n+1}\Big)'-t\Big(\sum_{n=0}^{\infty}d^{(r)}_n(x)t^n\Big)'
\\&\q-(x+r+1)\sum_{n=0}^{\infty}d^{(r)}_n(x)t^n-(x-r)\sum_{n=0}^{\infty}d^{(r)}_n(-x)(-t)^n
\\&=(1-t)\Big(\sum_{n=0}^{\infty}d^{(r)}_n(x)t^n\Big)'-(x+r+1)(1+t)^{x-r}(1-t)^{-x-r-1}
\\&\q-(x-r)(1-t)^{-x-r}(1+t)^{x-r-1}
\\&=(1-t)((1+t)^{x-r}(1-t)^{-x-r-1})'-(x+r+1)(1+t)^{x-r}(1-t)^{-x-r-1}
\\&\q-(x-r)(1-t)^{-x-r}(1+t)^{x-r-1}
\\&=0.
\endalign$$
This yields the result.
 \pro{Corollary 3.1} For $n=0,1,2,\ldots$ we
have
$$(n+2r)d_{n-1}^{(r)}(x)=(n+r-x)d^{(r)}_n(x)+(-1)^n(x-r)d^{(r)}_n(-x).$$
\endpro
Proof. Combining (1.3) and Theorem 3.1 gives the result.
 \pro{Corollary 3.2} For $n=1,2,3,\ldots$ we have
$$d^{(r)}_n\Ls 12=\cases 2(-1)^{\f{n-1}2}\b{-r-\f 32}{\f{n-1}2}&\t{if $2\nmid
n$,}
\\\f{2n+2r+1}n(-1)^{\f n2-1}\b{-r-\f 32}{\f n2-1}&\t{if $2\mid n$.}
\endcases$$
\endpro
Proof. By [5, Corollary 2.1],
$$d_n^{(r)}\Big(-\f 12\Big)=\cases 0&\t{if $2\nmid n$,}
\\(-1)^{\f n2}\b{-r-\f 12}{\f n2}&\t{if $2\mid n$.}\endcases\tag 3.1$$
Hence, from Theorem 3.1 we have
$$\align (-1)^n\Big(-\f 12-r\Big)d^{(r)}_n\Ls 12&=(n+1)d^{(r)}_{n+1}\Big(-\f
12\Big)-\Big(r+n+\f 12\Big)d^{(r)}_n\Big(-\f 12\Big)
\\&=\cases (n+1)(-1)^{\f{n+1}2}\b{-r-\f 12}{\f{n+1}2}&\t{if $2\nmid
n$,}
\\-(r+n+\f 12)(-1)^{\f n2}\b{-r-\f12}{\f n2}&\t{if $2\mid n$.}
\endcases\endalign$$
This yields the result for $r\not=-\f 12$. For $r=-\f 12$ from (1.2)
we have
$$\sum_{n=0}^{\infty}d_n^{(-\f
12)}\Ls 12t^n=\f{1+t}{1-t}=2+\f{2t}{1-t}=1+2\sum_{n=1}^{\infty}t^n\q
(|t|<1).$$ Thus, $d_n^{(-\f 12)}\sls 12=2$ for $n\ge 1$. This yields
the result for $r=-\f 12$. The proof is now complete.

 \pro{Corollary
3.3} For $r\not=-1$ and $n=0,1,2,\ldots$ we have
$$d^{(r)}_n(1)=\f{2n+1+(2-(-1)^n)r}{r+1}\b{r+[\f n2]}{[\f
n2]}.$$
\endpro
 Proof. By [5, (2.3) and
Corollary 2.2],
$$d_n^{(r)}(-1)=(-1)^nd_n^{(r)}(0)=(-1)^n\b{r+[\f n2]}{[\f n2]}.\tag
3.2$$ Taking $x=-1$ in Corollary 3.1 and then applying (3.2) gives
$$\align
(-1)^n(-1-r)d^{(r)}_n(1)&=(n+2r)d_{n-1}^{(r)}(-1)-(n+r+1)d_n^{(r)}(-1)
\\&=(n+2r)(-1)^{n-1}\b{r+[\f{n-1}2]}{[\f{n-1}2]}
-(n+r+1)(-1)^n\b{r+[\f n2]}{[\f n2]},\endalign$$ which yields the
result.

\pro{Theorem 3.2} For $n=1,2,3,\ldots$ we have
$$2d^{(r+\f 12)}_{n-1}\big(x-\f
12\big)=d^{(r)}_n(x)-(-1)^nd^{(r)}_n(-x)\tag 3.3$$ and
$$(n+1)d^{(r-\f 12)}_{n+1}\Big(x-\f
12\Big)=(x+r)d_n^{(r)}(x)+(x-r)(-1)^nd^{(r)}_n(-x).\tag 3.4$$
\endpro
Proof. Using (1.2) we see that for $|t|<1$,
$$\align
&\sum_{n=0}^{\infty}\big(d^{(r)}_n(x)-(-1)^nd^{(r)}_n(-x)\big)t^n
\\&=\f{(1+t)^{x-r}}{(1-t)^{x+r+1}}-\f{(1-t)^{-x-r}}{(1+t)^{-x+r+1}}
=\f{2t(1+t)^{x-r-1}}{(1-t)^{x+r+1}}=2t \f{(1+t)^{x-\f 12-(r+\f 12)}}
{(1-t)^{x-\f 12+r+\f 12+1}}
\\&=2\sum_{n=0}^{\infty}d^{(r+\f 12)}_n\big(x-\f 12\big)t^{n+1}.
\endalign$$
Comparing the coefficients of $t^n$ gives (3.3).
 Similarly,
using  (1.2) we see that for $|t|<1$,
$$\align
&\sum_{n=0}^{\infty}\big((n+1)d^{(r-\f 12)}_{n+1}\big(x-\f
12\big)-(x+r)d_n^{(r)}(x)-(x-r)(-1)^nd^{(r)}_n(-x)\big)t^n
\\&=\Big(\sum_{n=0}^{\infty}d^{(r-\f 12)}_{n+1}\big(x-\f
12\big)t^{n+1}\Big)'-(x+r)\sum_{n=0}^{\infty}d_n^{(r)}(x)t^n-(x-r)
\sum_{n=0}^{\infty}d_n^{(r)}(-x)(-t)^n
\\&=\big((1+t)^{x-r}(1-t)^{-x-r}\big)'-(x+r)(1+t)^{x-r}(1-t)^{-x-r-1}
-(x-r)(1+t)^{x-r-1}(1-t)^{-x-r}
\\&=0.
\endalign$$
This yields (3.4). Hence the theorem is proved.

\pro{Corollary 3.4} For $n=1,2,3,\ldots$ we have
$$2(x-r)d_{n-1}^{(r+\f 12)}\Big(x-\f 12\Big)+(n+1)d_{n+1}^{(r-\f
12)}\Big(x-\f 12\Big)=2xd_n^{(r)}(x).$$
\endpro
Proof. This is immediate from (3.3) and (3.4).

 \pro{Corollary 3.5} For $r\not=-1$ and $n=1,2,3,\ldots$ we have
 $$\align d_n^{(r)}\Big(\f 32\Big)=&\f 13\Big((n+1)\f{4n+6+(2+(-1)^n)(2r-1)}{2[\f{n+1}2]}
 \\&\q-(2r-3)
 \f{4n-2+(2+(-1)^n)(2r+1)}{2r+3}\Big)\b{r+\f 12+[\f{n-1}2]}{[\f{n-1}2]}.\endalign$$
\endpro
Proof. Taking $x=\f 32$ in Corollary 3.4 and then applying Corollary
3.3 gives
$$\align 3d_n^{(r)}\Ls 32&=2\Big(\f 32-r\Big)d_{n-1}^{(r+\f 12)}(1)+(n+1)d_{n+1}^{(r-\f 12)}(1)
\\&=(3-2r)\f{2(n-1)+1+(2-(-1)^{n-1})(r+\f 12)}{r+\f 12+1}\b{r+\f 12+[\f{n-1}2]}{[\f{n-1}2]}
\\&\q+(n+1)\f{2(n+1)+1+(2-(-1)^{n+1})(r-\f 12)}{r-\f 12+1}
\b{r-\f 12+[\f{n+1}2]}{[\f{n+1}2]},\endalign$$ which yields the
result. \pro{Theorem 3.3} For $n=0,1,2,\ldots$ we have
$$(x-r-1)d_n^{(r)}(1-x)=(x+r)(-1)^nd_n^{(r)}(x)-(2n+2r+1)d_n^{(r)}(-x).$$
\endpro
Proof. By (1.2), for $|t|<1$,
$$\align
&\sum_{n=0}^{\infty}\big((x+r)(-1)^nd_n^{(r)}(x)-(2n+2r+1)d_n^{(r)}(-x)\big)t^n
\\&=(x+r)\sum_{n=0}^{\infty}d_n^{(r)}(x)(-t)^n-(2r+1)\sum_{n=0}^{\infty}d_n^{(r)}(-x)t^n
-2\sum_{n=0}^{\infty}nd_n^{(r)}(-x)t^n
\\&=(x+r)(1+t)^{-x-r-1}(1-t)^{x-r}-(2r+1)(1+t)^{-x-r}(1-t)^{x-r-1}
-2t\Big(\sum_{n=0}^{\infty}d_n^{(r)}(-x)t^n\Big)'
\\&=(x+r)(1+t)^{-x-r-1}(1-t)^{x-r}-(2r+1)(1+t)^{-x-r}(1-t)^{x-r-1}
\\&\q-2t\big((1+t)^{-x-r}(1-t)^{x-r-1}\big)'
\\&=(x+r)(1+t)^{-x-r-1}(1-t)^{x-r}-(2r+1)(1+t)^{-x-r}(1-t)^{x-r-1}
\\&\q+2(x+r)t(1+t)^{-x-r-1}(1-t)^{x-r-1}+2(x-r-1)t(1+t)^{-x-r}(1-t)^{x-r-2}
\\&=(1+t)^{-x-r-1}(1-t)^{x-r-2}\big((x+r)(1-t)^2-(2r+1)(1-t^2)
\\&\q+2(x+r)t(1-t)+2(x-r-1)t(1+t)\big)
\\&=(x-r-1)(1+t)^{-x-r+1}(1-t)^{x-r-2}
\\&=(x-r-1)\sum_{n=0}^{\infty} d_n^{(r)}(1-x)t^n.
\endalign$$
Now comparing the coefficients of $t^n$ on both sides gives the
result.

\pro{Corollary 3.6} For $r\not=-1,-2$ and $n=0,1,2,\ldots$ we have
$$d_n^{(r)}(2)=\f{(3-2(-1)^n)r^2+(4-(-1)^n)(2n+1)r+4n^2+4n+2}{(r+1)(r+2)}\b{r+[\f
n2]}{[\f n2]}.$$
\endpro
Proof. Putting $x=-1$ in Theorem 3.3 and then applying (3.2) and
Corollary 3.3 we see that
$$\align
-(r+2)d_n^{(r)}(2)&=(r-1)(-1)^nd_n^{(r)}(-1)-(2n+2r+1)d_n^{(r)}(1)
\\&=(r-1)\b{r+[\f n2]}{[\f n2]}-(2n+2r+1)
\f{2n+1+(2-(-1)^n)r}{r+1}\b{r+[\f n2]}{[\f n2]}.\endalign$$ This
yields the result.

\section*{4. New proof of (1.4) and a formula for $M_n(x;b,-1)^2$}
\par For two positive integers $r$ and $s$ let the hypergeometric series
$_rF_s$ be given by
$$_rF_s\left(\begin{array}{rrr}a_1,&\ldots,&a_r\\b_1,&\ldots,&b_s
\end{array}\bigg|z\right)
=\sum_{k=0}^\infty\f{(a_1)_k\cdots(a_r)_k}{(b_1)_k\cdots(b_s)_k}
\cdot\f{z^k}{k!}.\tag 4.1$$ Recall that
$\f{(a)_k}{k!}=(-1)^k\b{-a}k=\b{a+k-1}k$ for $k=0,1,2,\ldots$. By
(2.1), for $r\not\in\{-\f 12,-\f 22,-\f 32,\ldots\}$,
$$d_n^{(r)}(x)=\f{(2r+1)_n}{n!} \sum_{k=0}^n
\f{(-n)_k(r-x)_k}{(2r+1)_k\cdot k!}\cdot 2^k =\f{(2r+1)_n}{n!}\
_2F_1\left(\begin{array}{rr}-n,&r-x\\&2r+1\end{array}\biggm|2\right).\tag
4.2$$
Since $d_n^{(r)}(x)=(-1)^nd_n^{(r)}(-1-x)$ we also have
$$d_n^{(r)}(x)=(-1)^n\f{(2r+1)_n}{n!}\
_2F_1\left(\begin{array}{rr}-n,&r+1+x\\&2r+1\end{array}\biggm|2\right).$$
Hence
$$d_n^{(r)}(x)^2=(-1)^n\f{((2r+1)_n)^2}{n!^2}
\
_2F_1\left(\begin{array}{rr}-n,&r+1+x\\&2r+1\end{array}\biggm|2\right)
\
_2F_1\left(\begin{array}{rr}-n,&r-x\\&2r+1\end{array}\biggm|2\right).\tag
4.3$$ From [4, (2.5.32)] we know that for $|z|<1$,
$$_2F_1\left(\begin{array}{rr}a,&b\\&c\end{array}\bigg|z\right){}_2F_1
\left(\begin{array}{rr}a,&c-b\\&c\end{array}\bigg|z\right)=(1-z)^{-a}
{}_4F_3\left(\begin{array}{rrrr}a,&b,&c-a,&c-b\\&c,&\f
c2,&\f{c+1}2\end{array}\bigg|\f{z^2}{4(z-1)}\right).\tag 4.4$$ When
$-a$ is a positive integer, the above three hypergeometric series
reduce to finite sums. Hence (4.4) is also true. Now taking $a=-n,\
b=r+1+x,\ c=2r+1$ and $z=2$ in (4.4) and then applying (4.3) yields
$$ d_n^{(r)}(x)^2=\f{((2r+1)_n)^2}{n!^2}
\ _4F_3\left(\begin{array}{rrrr}-n,&r+1+x,&2r+1+n,&r-x\\&2r+1,&
\f{2r+1}2,&r+1\end{array}\biggm| 1\right).\tag 4.5$$ That is,
$$ d_n^{(r)}(x)^2=\b{-2r-1}n^2 \sum_{k=0}^n\f{\b
nk\b{-r-1-x}k\b{-2r-1-n}k\b{x-r}k} {\b{-2r-1}k
\b{-r-1/2}k\b{-r-1}k}(-1)^k.$$
 Since $$\align &\b{-2r-1}n\b
nk\b{-2r-1-n}k\\&=\b{-2r-1}k\b{-2r-1-k}{n-k}\b{-2r-1-n}k=
\b{-2r-1}k\b{-2r-1-k}k\b{-2r-1-2k}{n-k},\endalign$$ we get
$$\align d_n^{(r)}(x)^2&=\b{-2r-1}n
\sum_{k=0}^n\f{\b{-r-1-x}k\b{x-r}k\b{-2r-1-k}k\b{-2r-1-2k}{n-k}} {
\b{-r-1/2}k\b{-r-1}k}(-1)^k
\\&=\b{n+2r}n\sum_{k=0}^n\f{\b{x+r+k}k\b{x-r}k\b{-2r-1-k}k\b{n+2r+k}{n-k}}
{\b{-r-1/2}k\b{-r-1}k}(-1)^k.\endalign$$
Note that
$$\align \b{-r-\f12}k&=\f{-\f{2r+1}2(-\f{2r+3}2)\cdots(-\f{2r+2k-1}2)}{k!}
=\f1{(-2)^k}\cdot\f{(2r+1)(2r+2)\cdots(2r+2k)}{k!\cdot
2^k\cdot(r+1)(r+2)\cdots(r+k)}\\&=\f1{4^k\cdot
k!}\cdot\f{(-2r-1)(-2r-2)\cdots(-2r-2k)}{(-r-1)(-r-2)\cdots(-r-k)}\\
&=\f{\b{-2r-1}{2k}\b{2k}k}{4^k\b{-r-1}k}
=\f1{4^k\b{-r-1}k}\b{-2r-1}k\b{-2r-1-k}k\\
&=\f1{(-4)^k\b{-r-1}k}\b{2r+k}k\b{-2r-1-k}k.
\endalign
$$
We then obtain
$$d_n^{(r)}(x)^2=\b{n+2r}n\sum_{k=0}^n\f{\b{x-r}k\b{x+r+k}k\b{n+2r+k}{n-k}}
{\b{2r+k}k}4^k,\tag 4.6$$ as given in [5, Theorem 2.6].
\par By Corollary 2.4,
$$d_n^{(\f{b-1}2)}\big(x+\f{b-1}2\big)
=\f{(b)_n}{n!}M_n(x;b,-1)=\b{b+n-1}nM_n(x;b,-1).\tag 4.7$$ Combining
(4.6) and (4.7) gives for $b\notin\{0,-1,-2,\ldots\}$,
$$M_n(x;b,-1)^2=\f 1{\b{b+n-1}n^2}d_n^{(\f{b-1}2)}\big(x+\f{b-1}2\big)^2
=\f 1{\b{b+n-1}n}\sum_{k=0}^n\f{\b
xk\b{x+b-1+k}k\b{n+k+b-1}{n-k}}{\b{b-1+k}k}4^k.\tag 4.8$$ That is,
$$\Big(\sum_{k=0}^n\b nk\f{\b xk}{\b{b-1+k}k}2^k\Big)^2
=\f 1{\b{b+n-1}n}\sum_{k=0}^n\f{\b
xk\b{x+b-1+k}k\b{n+k+b-1}{n-k}}{\b{b-1+k}k}4^k.\tag 4.9$$ Note that
$\b{-a}k=(-1)^k\b{a+k-1}k$. Replacing $b$ with $-a$ in (4.9) gives
the following result.

\pro{Theorem 4.1} Let $a\not\in\{0,1,2,\ldots\}$ and
 $n\in\{0,1,2,\ldots\}$. Then
 $$\Big(\sum_{k=0}^n\f{\b nk\b xk}{\b ak}(-2)^k\Big)^2
 =\f {(-1)^n}{\b
 an}\sum_{k=0}^n\f{\b xk\b{a-x}k\b{n+k-a-1}{n-k}}{\b ak}4^k.$$
 \endpro
 \pro{Corollary 4.1} Let $a\not\in\{0,1,2,\ldots\}$ and
 $n\in\{0,1,2,\ldots\}$. Then
 $$\Big(\sum_{k=0}^n\f{\b nk}{\b ak}2^k\Big)^2
 =\f 1{\b
 an}\sum_{k=0}^n(-1)^{n-k}\f{a+1}{a+1-k}\b{n+k-a-1}{n-k}4^k.$$
 In particular, for $a=-\f 12$ we have
 $$\Big(\sum_{k=0}^n\f{\b nk}{\b {2k}k}(-8)^k\Big)^2
 =\f 1{\b{2n}n}\sum_{k=0}^n(-1)^k\f 1{1-2k}\b{n+k-\f 12}{n-k}4^{n+k}.$$
 \endpro
 Proof. Note that $\b{-1}k=(-1)^k$. Taking $x=-1$ in Theorem 4.1
 yields the first part. It is well known that $\b{-\f
 12}k=\b{2k}k(-4)^{-k}$. Thus the remaining part follows.
 \pro{Corollary 4.2} Let
 $n\in\{0,1,2,\ldots\}$. Then
 $$\Big(\sum_{k=0}^n\b nk\b xk\f{2^k}{k+1}\Big)^2
 =\f 1{n+1}\sum_{k=0}^n\b xk\b{-2-x}k\b{n+k+1}{2k+1}\f{(-4)^k}{k+1}.$$
 \endpro
 Proof. Since $\b{-2}k=(-1)^k(k+1)$, taking $a=-2$ in Theorem 4.1
 yields the result.

\section*{5. Inequalities involving $d_n^{(r)}(x)$}
\pro{Theorem 5.1} Let $r>-\f 12$, $x\not=-\f 12$ and
$n\in\{2,3,4,\ldots\}$. Then
$$\f{d_n^{(r)}(x)d_{n-1}^{(r)}(x)}{1+2x}\ge \f 1n\Big(\b{2r+n-1}{n-1}+d_{n-1}^{(r)}(x)^2\Big)>
0.$$
\endpro
Proof. Putting $y=x$ in [5, Theorem 2.8] gives
$$(1+2x)\sum_{k=0}^{n-1}
\f{(2r+k+1)\cdots(2r+n)}{(k+1)\cdots
n}d_k^{(r)}(x)^2=(n+2r)d_n^{(r)}(x)d_{n-1}^{(r)}(x).\tag 5.1$$ Note
that $d_0^{(r)}(x)=1$ and $\f{(2r+k+1)\cdots(2r+n)}{(k+1)\cdots
n}d_k^{(r)}(x)^2\ge 0$. The result follows from (5.1).
\pro{Corollary 5.1} Let $r>-\f 12$, $x<-\f 12$ and
$n\in\{0,1,2,\ldots\}$. Then $(-1)^nd_n^{(r)}(x)>0$.
\endpro
Proof. Since $d_0^{(r)}(x)=1$ and $d_1^{(r)}(x)=1+2x$, the result is
true for $n=0,1$. Now assume $n\ge 2$. From Theorem 5.1 we know that
$d_n^{(r)}(x)d_{n-1}^{(r)}(x)<0$. Hence the result follows.

\pro{Corollary 5.2} Let $r>-\f 12$, $x>-\f 12$ and
$n\in\{2,3,4,\ldots\}$. Then $d_n^{(r)}(x)>\f{(2x+1)^n}{n!}>0.$
\endpro
Proof.  Observe that $d_0^{(r)}(x)=1>0$ and $d_1^{(r)}(x)=1+2x>0$.
By Theorem 5.1, $d_n^{(r)}(x)d_{n-1}^{(r)}(x)>0$ for $n\ge 2$. Hence
$d_k^{(r)}(x)>0$ for $k=0,1,2,\ldots$. For $k\ge 2$, using Theorem
5.1 we see that
$$\f{d_k^{(r)}(x)d_{k-1}^{(r)}(x)}{1+2x}>\f{d_{k-1}^{(r)}(x)^2}k\qtq{and so}
\f{d_k^{(r)}(x)}{d_{k-1}^{(r)}(x)}>\f{1+2x}k.$$ Hence
$$d_n^{(r)}(x)=d_1^{(r)}(x)\prod_{k=2}^n\f{d_k^{(r)}(x)}{d_{k-1}^{(r)}(x)}
>(1+2x)\prod_{k=2}^n\f{1+2x}k=\f{(2x+1)^n}{n!}.$$
This proves the corollary.
\par\q
\par Based on calculations with Maple we pose the following conjecture.
 \pro{Conjecture 5.1} Let $r\ge 0$ and $-1\le x\le 0$. Then
$$(-1)^n(d_n^{(r)}(x)^2-d_{n+1}^{(r)}(x)d_{n-1}^{(r)}(x))>0\qtq{for}n=1,2,3,\ldots.$$
\endpro
\par\q
\newline{\bf Acknowledgments} The author is supported by the National Natural Science Foundation of China
 (Grant No. 11771173).

\end{document}